\documentclass[12pt,a4paper]{article}
\usepackage{amsmath}
\usepackage{amsthm}
\usepackage{amscd}
\usepackage{amssymb}
\usepackage{amsfonts}
\usepackage{color}
\usepackage{latexsym}
\usepackage{graphicx}
\newtheorem{theorem}{Theorem}
\newtheorem{lemma}[theorem]{Lemma}
\newtheorem{corollary}[theorem]{Corollary}
\newtheorem{proposition}[theorem]{Proposition}

\theoremstyle{definition}
\newtheorem{definition}[theorem]{Definition}
\newtheorem{example}[theorem]{Example}

\theoremstyle{remark}
\newtheorem*{remark}{Remark}
\newtheorem*{nota}{Notation}

\newtheorem*{ack}{Acknowledgements}

\numberwithin{equation}{section}
\newcounter{temp}

\makeatletter
\def\square{\RIfM@\bgroup\else$\bgroup\aftergroup$\fi
\vcenter{\hrule\hbox{\vrule\@height.6em\kern.6em\vrule}\hrule}\egroup}
\makeatother

\DeclareMathOperator{\Spec}{Spec}
\DeclareMathOperator{\Proj}{Proj}
\DeclareMathOperator{\Hom}{Hom}
\DeclareMathOperator{\Der}{Der}
\DeclareMathOperator{\Image}{Im}

\DeclareMathOperator{\cone}{cone}
\DeclareMathOperator{\st}{st}
\DeclareMathOperator{\link}{lk}

\reversemarginpar
\newcommand{\kkk}[1]{}
\newcommand{\ktrash}[1]{{}}
\definecolor{Klaus}{rgb}{0.99,0.05,0.83}

\newcommand{\klabel}[1]{\label{#1}\kkk{#1}}

\newcommand{\N}{\mathbb{N}}
\newcommand{\Z}{\mathbb{Z}}

\newcommand{\C}{\mathbb{C}}

\newcommand{\A}{\mathbb{A}}
\newcommand{\PP}{\mathbb{P}}

\newcommand{\kk}{\C}
\newcommand{\kx}{\mathbf{x}}
\newcommand{\kp}{\mathbf{p}}

\newcommand{\ka}{\mathbf{a}}
\newcommand{\kb}{\mathbf{b}}
\newcommand{\kc}{\mathbf{c}}
\newcommand{\MX}{\Delta_n\setminus X}
\newcommand{\innt}{\operatorname{int}}
\newcommand{\lambdaM}{\mu}
\newcommand{\lambdaN}{\lambda}
\newcommand{\surj}{\rightarrow\hspace{-1.2em}\rightarrow}

\newcommand{\kst}{\,|\;}
\newcommand{\kSt}{\,\big|\;}

\newcommand{\kss}{\scriptscriptstyle}
\newcommand{\kbb}{{\kss \bullet}}

\begin{document}

\title{Cotangent cohomology of Stanley-Reisner rings}
\author{Klaus Altmann \and Jan Arthur Christophersen}
\date{}
\maketitle
\begin{abstract}
Simplicial complexes $X$ provide
commutative rings $A(X)$ via the Stanley-Reisner construction.
We calculated the cotangent cohomology,
i.e., $T^1$ and $T^2$ of $A(X)$ in terms of $X$. These modules provide
information about the deformation theory of the algebro geometric objects
assigned to $X$.
\end{abstract}


\section{Introduction}\klabel{pre}
Denote $[n] := \{0,\ldots,n\}$ and let 
$\Delta_{n}:=2^{[n]}$ be the full simplex. 
A simplicial complex $X\subseteq \Delta_{n}$ with 
vertex set $[X]\subseteq [n]$
gives rise to an ideal
\[
I_X:=\langle\kx^p\kst p\in \Delta_{n}\setminus X\rangle
\subseteq \kk[x_0,\dots,x_n] =: P.
\]
The {\it Stanley-Reisner ring} is then $A_X=P/I_X$.  
We can associate the schemes ${\mathbb A}(X)=\Spec A_X$ and
${\mathbb P}(X) = \Proj A_X$ with these rings.
The latter looks like $X$ itself -- its simplices have just been replaced
by projective spaces.
\par

For each $\kk$-algebra $A$, there is a cohomology theory
providing modules $T^i_A$, see e.g.\
\cite{an:hom} or \cite{la:for}.
However, only $T^1_A$ and $T^2_A$ are relevant
for the deformation theory of $\Spec A$. 
The module $T^1_A$ collects the infinitesimal deformations of $A$
and  $T^2_A$ contains the obstructions for lifting these
deformations to decent parameter spaces. 
Eventually, $\Der_{\kk}(A,A)$ will be called $T^0_A$.
\par

The main result of the present paper is Theorem \ref{HNN};
it provides the modules $T^i$ ($i=1,2$)
for Stanley-Reisner rings $A_X$ in terms of the geometry of
the original simplicial complex $X$: 
The $T^i$ are $\Z^{[n]}$-graded, i.e., each degree $\kc$
corresponds to a monomial
of the quotient field of the ambient ring $P$. 
Splitting $\kc=\ka-\kb$ in its positive and negative part,
one obtains disjoint subsets $a,b\subseteq [n]$ as their respective supports.
They give rise to certain subsets
$N_{a-b}, \widetilde{N}_{a-b}\subseteq X$, 
cf.\ Section \ref{TComplex} (right before Lemma \ref{Imphi}),
and it is the cohomology
of their geometric realizations 
which provides the homogeneous
part $T^i_\kc$ of $T^i$.
%
\par

One has to be a little careful with the geometric realization
of a subset $N\subseteq X$ which is not necessarily a subcomplex;
in particular it depends on whether $\emptyset\in N$ or not
-- see Section \ref{TGeometry} for the definition.
In Theorem \ref{topopen}, we present a version of our $T^i$-formula that
uses only open subsets of $X$ or certain nice subcomplexes.
We have chosen a non-trivial example (introduced as
Example \ref{oct} in Section \ref{TComplex}) to illustrate the theory.
In particular, it is spread (in eight parts) throughout the text.
\par

In the case that $|X|$ is a homological sphere, Ishida and Oda have proven 
Theorem \ref{HNN} in \cite{iso:tor} using torus embeddings. 
Moreover, in \cite{sym:def}, Symms computes 
$\Hom(I_X,A_X)$
and $T^2_0$ when $|X|$ is a
$2$-dimensional manifold possibly with boundary.
Our method is straightforward and allows us to get the $T^i$ for all 
Stanley-Reisner algebras. We also compute the cup product $T^1\times T^1\to
T^2$ and the localization maps.
Theorem~\ref{inject} states that they are always injective.
\par

Information about the $T^i$, the cup product, and their behavior under
localization makes it possible to investigate the deformation theory 
of $\A(X)$ and $\PP(X)$.
In fact, this paper was originally motivated by a question from Sorin
Popescu about the smoothability of ${\mathbb P}(X)$ when $|X|\approx
S^n$ as this would have applications for degenerations of Calabi-Yau
manifolds.  
In the forthcoming paper \cite{ac:def}
we apply our results to the case when $X$ is a combinatorial manifold,
e.g., a sphere.
Here we can give very explicit results and a good understanding of the
deformations of $\A(X)$ and $\PP(X)$.
\par

\begin{ack} 
The first author is grateful for financial support from the University
of Oslo during his visits.  The second author is grateful for
financial support from the universities FU and HU Berlin during his visits.
\end{ack}


\section{Cotangent cohomology in terms of the complex}\klabel{TComplex}

\begin{nota}\klabel{exp}
We will often work in the polynomial ring $\kk[\kx]=\kk[x_0,\dots,x_n]$.
Monomials are written as $\kx^{\kp}\in\kk[\kx]$ with exponents 
$\kp\in \N^{n+1}$. The support of $\kp$ is defined as
$p:=\{i\in[n]\kst  \kp_i\neq 0\}$. On the other hand,
subsets $p\subseteq[n]$ will allways be identified
with their characteristic vector $\kp\in\{0,1\}^{n+1}$.
\end{nota}

Let $X\subseteq \Delta_{n}$ be a simplicial complex 
and denote by
$\{e_{p}\kst  p\in \MX\}$ a basis for 
$P^{\left| \MX\right|}$ parametrizing the generators
$\kx^p$ of $I_X$.
The generating relations among them are
\[
R_{p,q}:=\quad {\kx^{q\setminus p}}\,e_p-
{\kx^{p\setminus q}}\,e_q\, .
\]
The relations among these relations are
\[
R_{p,q,r}:\quad \kx^{r\setminus (p\cup q)}\,e_{p,q}-
\kx^{q\setminus (p\cup r)}\,e_{p,r}+
\kx^{p\setminus (q\cup r)}\,e_{q,r}\, .
\]


\begin{remark} What we have just described is a special case of the
so called Taylor resolution -- a construction of a free, but in general
not minimal, resolution of any monomial ideal. For a description
and proof of exactness see e.g.\ \cite{bps:mon}.
\end{remark}

\newcommand{\XD}{D}
\begin{example}\klabel{oct}
The following simplicial complex $\XD$ will serve as a running example
throughout the text: With vertex set $\{x,y,z,x',y',z'\}$,
we define $\XD\subseteq\Delta_5$ to be the union of the
octahedron with the 8 maximal faces
$(x^{(')}y^{(')}z^{(')})$ and the 3 diagonals
$(xx')$, $(yy')$, and $(zz')$. Hence, the set $\Delta_5\setminus\XD$ 
providing the 
generators of $I_\XD$ consists of all $p\subseteq\{x,y,z,x',y',z'\}$ with
$\#p\geq 3$ and containing at least one letter twice.
\end{example}

$$
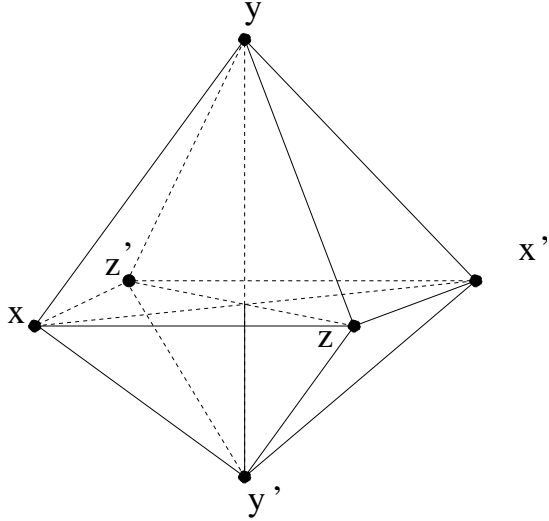
\begin{picture}(0,0)%
\includegraphics{diamant.pstex}%
\end{picture}%
\setlength{\unitlength}{1657sp}%
\begingroup\makeatletter\ifx\SetFigFont\undefined
\def\x#1#2#3#4#5#6#7\relax{\def\x{#1#2#3#4#5#6}}%
\expandafter\x\fmtname xxxxxx\relax \def\y{splain}%
\ifx\x\y   
\gdef\SetFigFont#1#2#3{%
  \ifnum #1<17\tiny\else \ifnum #1<20\small\else
  \ifnum #1<24\normalsize\else \ifnum #1<29\large\else
  \ifnum #1<34\Large\else \ifnum #1<41\LARGE\else
     \huge\fi\fi\fi\fi\fi\fi
  \csname #3\endcsname}%
\else
\gdef\SetFigFont#1#2#3{\begingroup
  \count@#1\relax \ifnum 25<\count@\count@25\fi
  \def\x{\endgroup\@setsize\SetFigFont{#2pt}}%
  \expandafter\x
    \csname \romannumeral\the\count@ pt\expandafter\endcsname
    \csname @\romannumeral\the\count@ pt\endcsname
  \csname #3\endcsname}%
\fi
\fi\endgroup
\begin{picture}(7967,7731)(2521,-7414)
\end{picture}
$$

In general, for a finitely generated $\kk$-algebra $A$,
the modules $T^i_A$ allow
the following ad hoc definitions:
Let $P=\kk[\kx]$ mapping onto $A$ so that $A\simeq P/I$
for an ideal $I$. Then $T^1_A$ is the cokernel of the natural map
$\Der_{\kk}(P,P) \to 
\Hom_P(I,A)$.
Moreover, if
\[
0 \to R \to P^m \stackrel{j}{\to} P \to A \to 0
\]
is an exact sequence presenting $A$ as a $P$ module and
$R_0:=\langle j(f) e- j(e) f \kst e,f\in P^m\rangle \subseteq R$
denotes the so-called Koszul relations, then
$R/{R_0}$ is an $A$ module and we obtain $T^2_A$ as the cokernel
of the induced map
$\,\Hom_P(P^m,A) \to \Hom_A(R/{R_0},A)$.

If $A=A_X$ is a Stanley-Reisner ring, then 
$A_X$ itself, its resolution, 
and all interesting $A_X$-modules 
such as the $T^i_A$ are $\Z^{n+1}$-graded;
just set
$\,\deg e_p=p$, $\deg R_{p,q}=p\cup q$, and
$\deg R_{p,q,r}=p\cup q\cup r$.
For an element  $\kc\in\Z^{n+1}$, we denote by 
\[
\Hom(I_X,A_X)_{\kc}
\hspace{1em}\mbox{and}\hspace{1em}
T^i_{\kc}(X):=T^i_{A_X,\kc}
\]
the homogeneous summands of the corresponding modules.
Let $\kc=\ka-\kb$ be the decomposition of $\kc$ in its
positive and negative part, i.e.,
$\,\ka,\kb\in\N^{n+1}$ with both elements having disjoint supports
$a$ and $b$, respectively.
This gives rise to the sets
\[
M_{a-b} :=\{p\in (\MX) \kst  (p\cup a)\setminus b \in X\}
\vspace{-1ex}
\]
and
\[
M^{(2)}_{a-b} := \{(p,q)\in M_{a-b}\times M_{a-b}\kst  
(p\cup q\cup a)\setminus b \in X\}\,.
\]
\par

{{\bf Example \ref{oct}.2} (continued) }
Let $p:=(xx'yy')$ and $q:=(xyy'z)$. Then $p,q\in M_{\emptyset-(yy')}$,
but $(p,q)\notin M^{(2)}_{\emptyset-(yy')}$.
Moreover, $(xx'yz), (xyz)\notin M_{\emptyset-(yy')}$, but
for different reasons.%
\vspace{1ex}
\par


\begin{lemma}\klabel{TiM}
Let $\,\kc=\ka-\kb$ as before.
The modules $\Hom(I_X,A_X)_{\kc}$ and $T^2_{\kc}(X)$ vanish unless
$\kb\in\{0,1\}^{n+1}$, i.e., $\kb=b$. Assuming $\kb=b$, these modules only
depend on the supports $a,b$.  
\vspace{0.5ex}\\
{\rm (i)} $\Hom(I_X,A_X)_{\kc}= \{\lambdaM:M_{a-b}\to\kk\kst 
\begin{array}[t]{@{}l}
\lambdaM(p)=0 \mbox{\rm\ if } b\not\subseteq p,\\
\lambdaM(p)=\lambdaM(q) \mbox{\rm\ if } (p,q)\in M^{(2)}_{a-b}\}\,.
\end{array}$
\vspace{0.5ex}\\
{\rm (ii)}
Elements of $\Hom(I_X,A_X)_{\kc}$ yield trivial deformations, i.e.,
belong to the image of $\Der_{\kk}(P,P)_{\kc}$, 
iff $\#(b)=1$ and
$\lambdaM(p)$ is a constant function.
\vspace{0.5ex}\\
{\rm (iii)}
$T^2_{\kc}(X)$ is the factor of
\vspace{-1.3ex}
\[
\{\lambdaM:M^{(2)}_{a-b}\to\kk\kst 
\begin{array}[t]{@{}l@{}}
\lambdaM \mbox{\rm\ is antisymmetric,}\;
\lambdaM(p,q)=0 
\mbox{\rm\ if } (p\cap q)\cup((p\cup q \cup a)\setminus b) \in X\\
\hspace{-2.5em}\mbox{\rm or if } b\not\subseteq p\cup q,\;
\lambdaM(p,q)-\lambdaM(p,r)+\lambdaM(q,r)=0 \mbox{\rm\ if } 
(p\cup q\cup r \cup a)\setminus b \in X\}
\vspace{-1.1ex}
\end{array}
\]
by
the subspace of functions $\lambdaM(p,q)=\lambdaM'(p)-\lambdaM'(q)$
with $\lambdaM'(m)=0$ if $b\not\subseteq m$.
\end{lemma}

\begin{proof} 
An element $\varphi\in \Hom(I_X,A_X)_{\kc}$ maps the generating monomials
$\kx^p$ to some $\lambdaM(p)\,\kx^{p+\ka-\kb}$
with $\lambdaM(p)\in\kk$. The condition that $(p+\ka-\kb)\in\N^{n+1}$
yields that $\lambdaM(p)=0$ unless $\kb=b$ and $b\subseteq p$.
On the other hand, if $\kx^{p+\ka-\kb}\in I_X$, then 
$\varphi(\kx^p)=0$, and the value $\lambdaM(p)$ does not matter at all.
Hence, we may restrict the knowledge of $\lambdaM$ to
$M_{a-b}\subseteq(\MX)$. The linearity of $\varphi$ translates into
the last condition in (i).
Eventually, the trivial deformations are spanned 
by $\varphi=\partial/\partial x_i$.\\
{(iii) } One obtains the description of 
$\Hom_A(R/{R_0},A_X)_{\kc}$
with the same arguments.
We should only remark that it is the 
Koszul relations $\,\kx^{q}e_p - \kx^{p} e_q\in R_0$ that are  
responsible for the vanishing of $\lambdaM(p,q)$ in case of
$\,(p\cap q)\cup((p\cup q \cup a)\setminus b) \in X$.
Afterwards, to get $T^2$, one needs to divide out the canonical generators
$D_m\in \Hom_{\kk[\kx]}(\kk[\kx]^{\MX},A_X)$ ($m \in \MX$). 
They have degree $-m$, and applied to
$R_{p,q}$, they yield non-trivial values 
$\kx^{q\setminus p}$ and $-\kx^{p\setminus q}$ only if $p=m$ or $q=m$,
respectively. Hence, in degree $\kc$, the map
$\kx^{m+\ka-b}D_m$ yields $\lambdaM(p,q)=\lambdaM_m(p)-\lambdaM_m(q)$
with $\lambdaM_m$ denoting the characteristic function of $m$.
On the other hand, this contribution requires $b\subseteq m$.
\end{proof} 

Of course, we are building some sort of cohomology to
describe the graded pieces of $T^i(X)$. However, in the previous lemma, 
the Koszul condition does not seem to fit. 
Surprisingly, this problem will be overcome by performing
kind of a $\kc$-shift.
Let
\vspace{-0.5ex}
\[
N_{a-b}(X):=\{f\in X \kst  a\subseteq f,\, f\cap
b=\emptyset ,\, f\cup b \notin X\}
\vspace{-1.0ex}
\]
and
\vspace{-0.5ex}
\[
{\widetilde N}_{a-b}(X):= \{f\in N_{a-b} \kst  
\text{$\exists$
$b^\prime \subset b$ with $f\cup b^\prime \notin X$}\}\, .
\vspace{1ex}
\]
\par


\begin{lemma} \klabel{Imphi}
Let $\Phi:M_{a-b}\to X$ be the application
$\Phi(p)=(p\cup a)\setminus b$.
\vspace{0.5ex}\\
{\rm (i)}
$\Image \Phi = N_{a-b}$. 
Moreover, an element $f\in N_{a-b}$ has a pre-image that does not
contain $b$ if and only if $f \in {\widetilde N}_{a-b}$.
\vspace{0.5ex}\\
%
{\rm (ii)} 
If $\,f,g\in N_{a-b}$ and $f \cup g
\in X$, then $f \cup g \in N_{a-b}$.  
If, moreover, $g\in {\widetilde N}_{a-b}$, then we even obtain
$f \cup g \in  {\widetilde N}_{a-b}$.
\end{lemma}

\begin{proof} 
(i) If $p\in M_{a-b}$, then $\Phi(p)\cup b\supseteq p\notin X$,
hence $\Phi(M_{a-b})\subseteq N_{a-b}$.
On the other hand, let $f\in N_{a-b}$ with $f\cup b'\notin X$
for some $b'\subseteq b$. It follows that $f\cup b'\in M_{a-b}$
and $\Phi(f\cup b')=f$. 
With $b':=b$, this implies $N_{a-b}\subseteq\Phi(M_{a-b})$;
with $b'$ being a proper subset of $b$, we obtain the 
${\widetilde N}_{a-b}$ statement.
The claims in (ii) are obvious.
\vspace{-1ex}
\end{proof}

{{\bf Example \ref{oct}.3} (continued) }
Considering the previous $p=(xx'yy')$ and $q=(xyy'z)$ from 
$M_{\emptyset-(yy')}$, we obtain 
$\,\Phi(p)=(xx')\in {\widetilde N}_{\emptyset-(yy')}(\XD)$ 
(note that $\,\Phi(p)=\Phi(xx'y)$), but
$\,\Phi(q)=(xz)\in N_{\emptyset-(yy')}(\XD) \setminus 
{\widetilde N}_{\emptyset-(yy')}(\XD)$.
\vspace{1ex}
\par

The map $\Phi:M_{a-b}\surj N_{a-b}$ of the previous lemma can easily be
extended to pairs. That is, with
$
N_{a-b}^{(2)} := \{(f,g)\in N_{a-b}\times N_{a-b}\kst  f \cup g \in X\}\,,
$
we also have a surjective application 
$\Phi:M_{a-b}^{(2)}\surj N_{a-b}^{(2)}$.
\par


\begin{proposition}\klabel{TiN}
Let $\,\kc=\ka-\kb$ as before.
The modules $\Hom(I_X,A_X)_{\kc}$ and $T^2_{\kc}(X)$ vanish unless
$\kb\in\{0,1\}^{n+1}$, i.e., $\kb=b$. Assuming $\kb=b$, these modules only
depend on the supports $a,b$.
\vspace{0.5ex}\\
{\rm (i)} 
$\Hom(I_X,A_X)_{\kc}= \{\lambdaN:N_{a-b}\to\kk\kst 
\begin{array}[t]{@{}l}
\lambdaN(f)=0 \mbox{\rm\ if } f\in {\widetilde N}_{a-b},\\
\lambdaN(f)=\lambdaN(g) \mbox{\rm\ if } f\cup g\in X\}\,.
\end{array}$
\vspace{0.5ex}\\
{\rm (ii)}
Elements of $\Hom(I_X,A_X)_{\kc}$ yield trivial deformations, i.e.,
belong to the image of $\Der_{\kk}(P,P)_{\kc}$, 
iff $\#(b)=1$ and
$\lambdaN(f)$ is a constant function.
\vspace{0.5ex}\\
{\rm (iii)}
$T^2_{\kc}(X)$ is the factor of the vector space of
the antisymmetric maps $\lambdaN:N^{(2)}_{a-b}\to\kk$ such that
\vspace{-1.3ex}
\[
\lambdaN(f,g)=0 \mbox{ if } f,g\in {\widetilde N}_{a-b}
\hspace{0.5em}\mbox{and}\hspace{0.5em}
\lambdaN(f,g)-\lambdaN(f,h)+\lambdaN(g,h)=0 \mbox{ if } 
(f\cup g\cup h)\in X
\vspace{-1.1ex}
\]
by the subspace $\{\lambdaN(f)-\lambdaN(g)\}$ with
$\lambdaN=0$ on ${\widetilde N}_{a-b}$.
\end{proposition}

\begin{proof}
Denote, just for this proof, the spaces given by (i) and (iii) of the previous
proposition by $\Hom(N)$ and $T^2(N)$, respectively.
Then we have to ascertain that 
pulling back via
$\Phi$ induces isomorphisms
$\Phi^\ast:\Hom(N)\stackrel{\sim}{\to}\Hom(M)$ and 
$\Phi^\ast:T^2(N)\stackrel{\sim}{\to} T^2(M)$
with $\Hom(M)$ and $T^2(M)$ being the corresponding spaces from
Lemma \ref{TiM}.
\vspace{0.5ex}\\
{\em Step 1. }
{\em The maps $\Phi^\ast$ are correctly defined: }\\
This is clear for the Hom case. For $T^2$, we set
$\lambdaM(p,q):=\lambdaN(\Phi p, \Phi q)$, and the only non-trivial task 
is to check the two
conditions that should lead to the vanishing of $\lambdaM(p,q)$. 
If $b\not\subseteq p\cup q$, then both 
$b\not\subseteq p$ and $b\not\subseteq q$, hence 
$\Phi(p),\Phi(q)\in {\widetilde N}$ by Lemma~\ref{Imphi}(i),
hence $\lambdaN(\Phi p, \Phi q)=0$.
On the other hand,
if $(p\cap q)\cup((p\cup q \cup a)\setminus b) \in X$, then 
we also obtain $b\not\subseteq p$ and $b\not\subseteq q$:
Otherwise, if $b\subseteq p$,
it would follow that $q\cap b \subseteq p\cap q$
and $q\setminus b\subseteq (p\cup q \cup a)\setminus b$,
thus, $q=(q\cap b)\cup (q\setminus b) \subseteq
(p\cap q)\cup((p\cup q \cup a)\setminus b)\in X$, 
but this contradicts $q\in M$.
\vspace{0.5ex}\\
{\em Step 2. }
{\em The maps $\Phi^\ast$ are injective: }\\
The Hom case follows from the surjectivity of $\Phi$.
For $T^2$, assume that 
$\lambdaN(\Phi p, \Phi q)=\lambdaM(p,q)=\lambdaM(p)-\lambdaM(q)$.
In particular, if $p,q$ belong to a
common fiber $\Phi^{-1}(f)$, then $\lambdaM(p)-\lambdaM(q)=0$.
Hence, $\lambdaM(p), \lambdaM(q)$ only depend on 
$\Phi(p)$ and $\Phi(q)$.
\vspace{0.5ex}\\
{\em Step 3. }
{\em The maps $\Phi^\ast$ are surjective: }\\
Let $\{\lambdaM(p)\}$ represent an element of $\Hom(M)$. Then, the property
that $\lambdaM(p)=\lambdaM(q)$ for $(p,q)\in M^{(2)}$ implies that $\mu(p)$
only depends on $\Phi(p)$. In particular, 
$\{\lambdaM(p)\}\in 
\Phi^\ast(\Hom(N))$.
\vspace{0.5ex}\\
To check the $T^2$ case,
we would like to proceed similarily with elements 
$\{\lambdaM(p,q)\}\in T^2(M)$. 
However, this requires a correction by coboundaries:
By the cocycle property of the $\lambdaM(p,q)$'s, 
we have to find $\{\lambdaM(p)\}$ such that
$\tilde{\lambdaM}(p,q):= 
\lambdaM(p,q) + (\lambdaM(p)-\lambdaM(q))$
vanishes if $p,q$
belong to a common fiber $\Phi^{-1}(f)$.
Using the cocycle property again, we see that
$\lambdaM(p):=\lambdaM(m_f,p)$ will almost do the job for any
fixed $m_f\in\Phi^{-1}(f)$; but we also have to ensure that 
$\lambdaM(p)=0$ whenever $b\not\subseteq p$.
This is done by proving the following\\
%
{\em Claim: } 
{\em Let $m,p\in M$ with $b\not\subseteq m$, $b \not\subseteq p$,
and $\Phi(m)\supseteq\Phi(p)$. Then $\lambdaM(m,p)=0$.}\\
With $f:=\Phi(m)$, we have that
$(m\cap p)\cup f= 
(m\cap p)\cup((m\cup p \cup a)\setminus b)$.
Thus, $(m\cap p)\cup f \notin X$ (or the Koszul condition 
in Lemma \ref{TiM}(iii) immediatly implies $\lambdaM(m,p)=0$).
In particular, $(m\cap p)\cup f\in  \Phi^{-1}(f)$, and
since 
$b\not\subseteq m\cup [(m\cap p)\cup f]$ and
$b\not\subseteq [(m\cap p)\cup f]\cup p$,
we obtain that
\vspace{-1ex}
\[
\lambdaM(m,p)\;=\;
\lambdaM(m,\,(m\cap p)\cup f) + \lambdaM((m\cap p)\cup f,\,p)
\;=\;0+0\;=\;0\,.
\vspace{-0.5ex}
\]
Eventually, it is possible to define
$\lambdaN(\Phi p, \Phi q):=\tilde{\lambdaM}(p,q)
=\lambdaM(m_{\Phi(p)}, m_{\Phi(q)})$, and it remains
to show its vanishing for 
$\Phi(p),\Phi(q)\in {\widetilde N}$.
Since,  by Lemma \ref{Imphi}(ii),
$\Phi(p)\cup \Phi(q)\in {\widetilde N}$,
we may assume that $\Phi(p)\subseteq \Phi(q)$. 
Now everything follows from applying the previous claim again.
\end{proof}


\begin{definition}\klabel{Uprop} 
A subset $Y\subseteq X$ of a simplicial complex $X$
has property U, or is a U subset, if 
\vspace{-1ex}
\begin{equation*}f,g\in Y \text{
and } f\cup g\in X \Rightarrow f\cup g\in Y\, .
\vspace{1ex}
\end{equation*}
\end{definition}

If $Y$ has this property, then we define the sets
\[
Y^{(k)}:=\{(f_0,\dots,f_k) \in Y^{k+1}\kst  f_0\cup\dots\cup f_k \in Y\}
\]
and the complex of $\kk$-vector spaces
\[
K^k(Y):=\{\lambda : Y^{(k)} \to \kk \kst  \lambda\text{ is
alternating}\}\subseteq \Lambda^{k+1}\big(\kk^Y\big)
\]
with the usual differential $d:K^{k-1}(Y)\to K^{k}(Y)$ defined by
\[
d(\lambda)(f_0,\dots,f_k):=\sum_{v=0}^k (-1)^v
\lambda(f_0,\dots,\hat{f_v},\dots,f_k)\,.
\]
By Lemma~\ref{Imphi}(ii), both $N_{a-b}$ and ${\widetilde N}_{a-b}$ are
U subsets.
Moreover,
there is a canonical surjection of complexes
$K^\bullet(N_{a-b})\twoheadrightarrow
K^\bullet({\widetilde N}_{a-b})$ leading to


\begin{corollary}\klabel{complex}
Assume $\mathbf{c}=\mathbf{a}-b$ with disjoint $a,b\in X$. 
Then
\begin{align*}
\Hom(I_X,A_X)_{\mathbf{c}}&\simeq H^0(\ker
(K^\bullet(N_{a-b}) \twoheadrightarrow
K^\bullet({\widetilde N}_{a-b})))\\
T^2_{\mathbf{c}}(X)&\simeq H^1(\ker
(K^\bullet(N_{a-b}) \twoheadrightarrow
K^\bullet({\widetilde N}_{a-b})))\,,
\end{align*}
and the trivial deformations inside $\Hom(I_X,A_X)_\kc$, i.e., those yielding
$0$ in $T^1_\kc(X)$, form a one-dimensional subspace whenever
$\#(b)=1$ (and are absent otherwise).
\end{corollary}

{{\bf Example \ref{oct}.4} (continued) }
We still consider the degree
$\,a=\emptyset$, $\,b=\{y,y'\}$ for $\XD$ being the octahedron with diagonals.
The set $N_{\emptyset-(yy')}(\XD)$ consists of the 4 vertices
$x^{(')}$, $z^{(')}$ and all the 6 edges connecting them.
However, only the interior of the edges $xx'$ and $zz'$ survive
in ${\widetilde N}_{\emptyset-(yy')}(\XD)$. 
\\

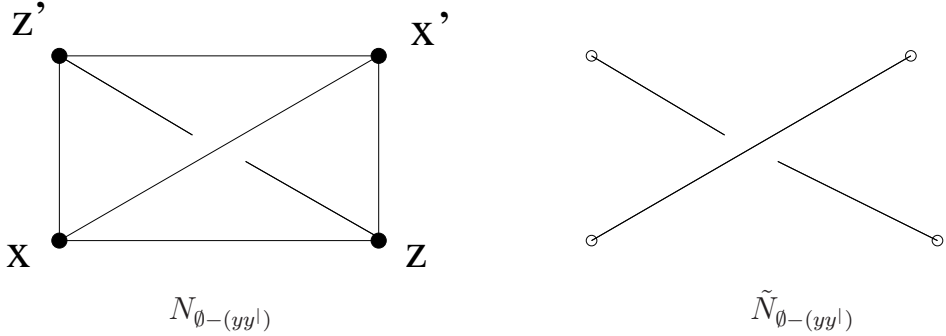
\begin{picture}(0,0)%
\includegraphics{two.pstex}%
\end{picture}%
\setlength{\unitlength}{2901sp}%
\begingroup\makeatletter\ifx\SetFigFont\undefined
\def\x#1#2#3#4#5#6#7\relax{\def\x{#1#2#3#4#5#6}}%
\expandafter\x\fmtname xxxxxx\relax \def\y{splain}%
\ifx\x\y   
\gdef\SetFigFont#1#2#3{%
  \ifnum #1<17\tiny\else \ifnum #1<20\small\else
  \ifnum #1<24\normalsize\else \ifnum #1<29\large\else
  \ifnum #1<34\Large\else \ifnum #1<41\LARGE\else
     \huge\fi\fi\fi\fi\fi\fi
  \csname #3\endcsname}%
\else
\gdef\SetFigFont#1#2#3{\begingroup
  \count@#1\relax \ifnum 25<\count@\count@25\fi
  \def\x{\endgroup\@setsize\SetFigFont{#2pt}}%
  \expandafter\x
    \csname \romannumeral\the\count@ pt\expandafter\endcsname
    \csname @\romannumeral\the\count@ pt\endcsname
  \csname #3\endcsname}%
\fi
\fi\endgroup
\begin{picture}(7928,2265)(1351,-2311)
\end{picture}

\hspace*{2cm} $N_{\emptyset - (yy^|)}$ \hspace*{6cm}  $\tilde{N}_{\emptyset - (yy^|)}$ 

To obtain elements of $\Hom(I_\XD,A_\XD)_{\mathbf{c}}$, we have to consider
maps $\lambda:N\to\kk$, i.e., each of the 4 vertices and 6 edges will be
assigned a value. The two conditions encoded
by ``$H^0$'' and ``$\ker$'' in the previous corollary mean
that $\lambda$ has to be both constant along the graph and zero
on $\innt xx'$ and $\innt zz'$. Hence, $\Hom(I_\XD,A_\XD)_{\mathbf{c}}=0$.
\par


\begin{remark}\klabel{cohom}
As we mentioned in the beginning, there is a general
cohomological definition of the $T^i$. 
Hence, it is no surprise that we
ended up with a cohomological description of these spaces in terms of $X$, too.
Moreover,
it should even be a challange to find a {\em direct} way 
to obtain the previous result (without touching elements).
If this involved a description of the so-called cotangent complex, one would
obtain important information about the deformation theory of both
$\A(X)$ and $\PP(X)$.
\end{remark}
\par


\section{Cotangent cohomology and the geometry of $X$}\klabel{TGeometry}

In the following, we will relate the previous description
of $T^i(X)$ with the geometry of the complex.
Let us start with some notation.
%
%
For $g \subseteq [n]$, denote by $\bar{g}:=2^g$
and $\partial g := \bar{g}\setminus \{g\}$
the full simplex and its boundary, respectively. 
The {\it join} $X\ast Y$ of two complexes $X$ and $Y$
is the complex defined by
\vspace{-1ex}
\[
X\ast Y := \{f\vee g : f\in X,\, g\in Y\}
\]
where $\vee$ means the disjoint union. If $f\in X$ is a face, we may
define
\vspace{-1ex}
\begin{itemize}
\item the {\it link} of $f$ in $X$; $\;\link(f,X):=\{g\in X: g\cap f =
\emptyset \text{ and } g\cup f\in X\}$,
\vspace{-1ex}
\item the {\it open star} of $f$ in $X$; 
$\;\st(f,X):=\{g\in X: f\subseteq g\}$, and
\vspace{-1ex}
\item the {\it closed star} of $f$ in $X$; 
$\;\overline{\st}(f,X):=\{g\in X: g\cup f\in X\}$.
\vspace{-1ex}
\end{itemize}
Notice that the closed star is the subcomplex $\overline{\st}(f,X) =
\bar{f}\ast \link(f,X)$.
%
Recall that the {\em geometric realization} of $X$, denoted $|X|$,
may be described by
\[
|X| = \big\{\alpha: [n] \to [0,1] \kSt \{i\kst \alpha(i) \ne 0\}\in X 
\mbox{ and  $\,\sum_i \alpha(i) = 1$} \big\}\, .
\]
To every non-empty $f\in X$, one assigns the {\em relatively open}
simplex $\langle f\rangle \subseteq |X|$;
\[
\langle f\rangle  = \{\alpha\in |X|\kst \alpha(i) \ne 0 \text{ if and only if }
i\in f \}\, .
\]
%
On the other hand, each subset $Y \subseteq X$ determines a
topological space
\[
\langle Y\rangle:=
\begin{cases}
        \bigcup_{f\in Y}\langle f\rangle& \text{if $\emptyset\not\in Y$}, \\
        \cone \left(\bigcup_{f\in Y}\langle f\rangle\right)& \text{if
        $\emptyset\in Y$}\, .
\end{cases}
\]
In particular, $\langle X\setminus\{\emptyset\}\rangle = |X|$ and
$\langle X\rangle = |\cone(X)|$ where $\cone(X)$ is the simplicial
complex $\Delta_0\ast X$.

Any subset $Y$ of $X$ is a poset with respect to inclusion and we may
construct the associated (normalized) order complex $Y^\prime$:
The vertices of $Y^\prime$ are the elements of $Y$ and the $k$-faces
of $Y^\prime$ are flags $f_{0}\subset f_{1}\subset\dots\subset f_{k}$
of $Y$-elements. If $Y$ is a complex, $Y^\prime$ is the barycentric
subdivision of $\cone(Y)$.

A complex and its barycentric subdivision have the same geometric
realization, so if $Y$ is a subcomplex of $X$, we have
$|Y^\prime|=\langle Y\rangle$.  This identity is obtained by sending a
vertex $f$ of $Y^\prime$ to the barycenter of $f$ in $\langle
f\rangle$ if $f\ne\emptyset$ and the vertex corresponding to
$\emptyset\in Y^\prime$ to the vertex of the cone.  For a general
subset $Y\subseteq X$, we only know that $|Y^\prime| \subseteq \langle
Y\rangle$ inside $|X^\prime|=\langle X\rangle$.

\begin{lemma}\klabel{retract}
If $Y\subseteq X$, then $|Y^\prime|$ is a deformation retract
of $\langle Y\rangle$. In particular, both sets have the same cohomology.
\end{lemma}

\begin{proof} If $f\in \cone(X)$, then we may identify $\langle f\rangle
\subset \langle X\rangle$ with the union of all $\langle F\rangle$ in
$|X^\prime|$ where $F=(f_{0}\subset f_{1}\subset\dots\subset f_{k})$
and $f_{k}=f$. Thus, $\langle Y\rangle$ as subset of $|X^\prime|$ is
the union of such $\langle F\rangle$ with $f_{k}\in Y$.  For such an
$F$ let $F_Y\leq F$ be the maximal subflag consisting only of faces
in $Y$.  Now we can continuously retract $\langle F\rangle\cup \langle
F_Y\rangle$ onto $\langle F_Y\rangle$ and this can be done
simultaneously for all $F$ belonging to the above union.
\end{proof}

For a subset $Y\subseteq X$, let $C^\bullet(Y)$ be the cochain complex of
$Y^\prime$.  We have an obvious inclusion of the $k$-flags in $Y$
into $Y^{(k)}$ given by $f_{0}\subset f_{1}\subset\dots\subset f_{k}
\mapsto (f_{0},f_{1},\dots ,f_{k})$.  This induces a surjection of
complexes $K^\bullet(Y) \rightarrow C^{\bullet }(Y)$ when $Y$
has property U.

\begin{lemma}\klabel{quasi}
The surjection $K^\bullet(Y) \rightarrow C^{\bullet }(Y)$ is a
quasi-isomorphism, 
i.e., it
induces an isomorphism in cohomology.
\end{lemma}

\begin{proof} We will prove the dual statement in homology using the
method of acyclic models (see e.g.\  \cite[4.2]{spa:alg}).  If $Y$ has
property U,
consider the simplicial complex $Y^{*}$ where the vertices are the
elements of $Y$ and a set of vertices $\{f_{0},\dots ,f_{k}\}$ is a
face if $f_{0}\cup\dots\cup f_{k}\in Y$.  If $C_\bullet(Y^{*})$ is the
chain complex of $Y^{*}$, then our $K^\bullet(Y)$ is the dual of
$C_{\bullet}(Y^{*})\otimes \kk$, and 
we are finished if we can 
prove that
$C_\bullet(Y^{*})$ is chain equivalent to $C_\bullet(Y^{\prime})$.\\
To this end, consider the set $\mathcal Y$ of U subsets of $X$ as a
category with inclusions as morphisms.  Let the models in $\mathcal Y$
be $\mathcal M=\{\bar{f}\cap Y\kst f\in Y,\, Y\in \mathcal Y\}$.  Finally,
define the two functors from $\mathcal Y$ to chain complexes by
$F^\prime(Y)=C_\bullet(Y^{\prime})$ and $F^*(Y)=C_\bullet(Y^{*})$.  We
must show that both $F^\prime$ and $F^*$ are free and acyclic with
respect to these models.
Now a basis element of $F^\prime(Y)$, $(f_{0}\subset
f_{1}\subset\dots\subset f_k)$, comes from $F^\prime(
\overline{f_{k}}\cap Y)$ and a basis element $\{f_{0},\dots ,f_{k}\}$
of $F^*(Y)$ comes from $F^*(\overline{f_{0}\cup\dots\cup f_{k}}\cap
Y)$, so both functors are free.  If $f\in Y$, one may check that
$(\bar{f}\cap Y)^{*}$ is a simplex and $(\bar{f}\cap Y)^\prime$ is a
cone over the vertex $\{f\}$.  Thus, in both cases the chain complexes
are acyclic.
\end{proof}

We may apply Lemma~\ref{retract} and \ref{quasi}
to $N_{a-b}$ and ${\widetilde N}_{a-b}$. 
Via the 5-Lemma, the result of Corollary~\ref{complex} translates into

\begin{theorem}\klabel{HNN}
Assume $\kc=\ka-\kb$ with $\ka,\kb\in\N^{n+1}$ having disjoint supports
$a,b\subseteq[n]$.
The homogeneous pieces in degree $\kc$ of the cotangent cohomology 
of the Stanley-Reisner ring $A_X$ vanish unless
$\kb\in\{0,1\}^{n+1}$, i.e., $\kb=b$. 
Assuming $\kb=b$, these modules only depend on the supports $a,b$,
and we have isomorphisms
\[
T^i_\kc(X)\;\simeq\; H^{i-1}(\langle
N_{a-b}\rangle,\langle{\widetilde N}_{a-b}\rangle,\kk)
\;\text{ for } i=1,2
\]
unless $b$ consists of a single vertex.  If $b$ consists of only one
vertex then the above formulae become true if we use the reduced
cohomology instead.\\
Moreover, replacing $T^1_\kc(X)$ by $\Hom(I_X,A_X)_\kc$ creates true
formulae without any need to distinguish between several cases.
\end{theorem}
\par

{{\bf Example \ref{oct}.5} (continued) }
In the previous session of Example \ref{oct}, we have seen that
the pair $(N_{\emptyset-(yy')}, \,{\widetilde N}_{\emptyset-(yy')})$
equals 
(complete graph on 4 vertices, 2 opposite edges).
Now, we may use that this is homotopy equivalent to
$(\unitlength=0.40pt
\hspace{0.2em}
\begin{picture}(150.00,20.00)(70.00,740.00)
\put(80.00,759.00){\line(1,0){120.00}}
\put(80.00,753.00){\line(1,0){120.00}}
\put(200.00,747.00){\line(-1,0){120.00}}
\put(200.00,741.00){\line(-1,0){120.00}}
\put(70.00,750.00){\circle*{15}}
\put(212.00,750.00){\circle*{15}}
\end{picture}
\hspace{0.2em},\hspace{0.4em}
\begin{picture}(30.00,20.00)(20.00,740.00)
\put(20.00,750.00){\circle*{15}}
\put(40.00,750.00){\circle*{15}}
\end{picture})$,
i.e.,
$T^1_\kc(\XD)=0$, but $\,\dim_\kk T^2_\kc(\XD)=4$.
\vspace{1ex}
\par

The previous theorem also provides information about
$T^0(X):=\Der_\kk(A_X,A_X)$. While this module has already been
described in \cite{bs:mod}, we would like to demonstrate its relation
to our techniques.
\par


\begin{corollary}\klabel{T0}
$\,T^0(X)=\bigoplus_{v=0}^n \mathfrak{a}_v\,\partial /\partial x_v$
where $\mathfrak{a}_v$ is the ideal of $A_X$ 
\vspace{0.3ex}
generated by the monomials $x^a$ with 
$\,\overline{\st}(a,X) \subseteq \overline{\st}(v,X)$.\\
In particular, $T^0(X)$ is generated, as a module, by
$\delta_v: = x_v \, \partial /\partial x_v$
if and only if 
every non-maximal $a\in X$ is properly contained in at least two
different faces.
\end{corollary}

\begin{proof}
$\,T^0(X)$ is the kernel of 
$\,\Der_\kk(\kk[\kx],A_X)\to \Hom(I_X,A_X)$.
Hence, since ${\widetilde N}_{a-\{v\}}=\emptyset$,
the previous theorem implies that
an element $\kx^\ka\,\partial/\partial x_v$ belongs to $T^0(X)$
if and only if 
$H^0(\langle N_{a-\{v\}}\rangle,\kk)=0$, 
i.e., iff $N_{a-\{v\}}=\emptyset$.
On the other hand, this means that, for every $f\in X$,
the conditions $a\subseteq f$ and $v\notin f$ imply $f\cup v\in X$.
Since the assumption $v\notin f$ can be omitted, this translates
into $\,\st(a,X)\subseteq \overline{\st}(v,X)$.\\
Finally, $T^0(X)$ is generated, as a module, by
$\delta_v = x_v \, \partial /\partial x_v$
if and only if
$\,\overline{\st}(a,X) \subseteq \overline{\st}(v,X)$
cannot happen for faces $a$ with $v\notin a$.
But this is equivalent to the condition formulated in the corollary.
\end{proof}


\section{Reduction to the $a=\emptyset$ case and localization}
\klabel{empty}

The set $N_{a-b}$ is empty unless $b\ne\emptyset$ and $a\in X$ is a
face.  Moreover, by the next proposition,
we may reduce all the calculations of $N_{a-b}$ and
$\widetilde{N}_{a-b}$, and therefore the $T^i$, to the case of
$a=\emptyset$ on a smaller complex.  
See Example \ref{En} for a demonstration of a consequent usage of this method.
\par


\begin{proposition}\klabel{aempty}
$\,T^i_\kc(X)=0$ for $\,i=1,2$ unless $\,a\in X$ and $\,b\subseteq
[\link(a)]$.  If $\,b\subseteq [\link(a)]$, then the map $f\mapsto
f\setminus a$ is a bijection
$N_{a-b}(X)\stackrel{\sim}{\to} N_{\emptyset -b}(\link(a))$ 
inducing isomorphisms
$T^i_{\emptyset -b}(\link(a)) \simeq T^i_{\ka-b}(X)$ for $i=1,2$.
\end{proposition}

\begin{proof} First assume that there is a vertex $v\in b
\setminus [\link(a)]$.  If $f\in N_{a-b}$, then $f\cup v\not\in X$
(otherwise $a\cup v\in X$ and $v\in \link(a)$). Thus,
$N_{a-b}={\widetilde N}_{a-b}$ unless $b=\{v\}$.  If $b=\{v\}$, then
$\widetilde{N}_{a-b}=\emptyset$ and $N_{a-b}=\st(a)$.  Thus, $\langle
N_{a-b}\rangle$ may be contracted to $\langle a\rangle$ and its
reduced cohomology is trivial, so $T^i_{a-b}(X)=0$ by
Theorem~\ref{HNN}.\\
It is a simple matter to check that the map between the $N$ sets is a
bijection (with inverse $g\mapsto g\cup a$) and that it restricts to a
bijection of the $\widetilde{N}$ subsets.  Since it clearly preserves
inclusions, it induces a simplicial isomorphism on the complexes
$N_{a-b}(X)^\prime \simeq N_{\emptyset -b}(\link(a))^\prime$. 
From Lemma~\ref{retract}, it follows
that we get an isomorphism in the
relative cohomology.
\vspace{-1ex}
\end{proof}

{{\bf Example \ref{oct}.6} (continued) }
Assume that $a=\{x\}$ for $\XD$ from Example \ref{oct}.
The link $\,\link(x)$ equals the boundary of the rectangle
$(yzy'z')$ plus the isolated point $x'$.
\vspace{0.3ex}\\
If $\#b\geq 3$, then there is always a proper subset
$b'\subset b$ with $b'\notin X$. In particular,
$N_{\emptyset-b}(\link(x))=\widetilde{N}_{\emptyset-b}(\link(x))$,
hence $T^1_{(x)-b}(\XD)=T^2_{(x)-b}(\XD)=0$.
\vspace{0.3ex}\\
The case $\#b=2$ does not provide any $T^i$ either, but if $b=\{x'\}$,
then $N_{\emptyset-(x')}(\link(x))$ is the boundary of the rectangle
and $\widetilde{N}_{\emptyset-(x')}(\link(x))=\emptyset$.
Hence, $T^1_{(x)-(x')}(\XD)=0$ and $\,\dim_\kk T^2_{(x)-(x')}(\XD)=1$.
Similarily, we obtain
$\,\dim_\kk T^1_{(x)-\ast}(\XD)=1$ and $T^2_{(x)-\ast}(\XD)=0$
where $\ast$ stands for any of the vertices $y^{(')}$, $z^{(')}$.
\vspace{1ex}
\par

The subsets $\langle N_{a-b} \rangle$ and
$\langle{\widetilde N}_{a-b}\rangle$ are in general neither open
nor closed in $|X|$.  In the case $a=\emptyset$ though, we may find open sets
retracting onto them.  
These sets are easier to define, but they
are not always easier to handle. However,
their openess often allows use of standard tools for calculating
cohomology.
Let 
\vspace{-0.5ex}
\[
U_{b} = U_{b}(X) :=\{f \in X\kst f\cup b \not\in X\}
\vspace{-0.5ex}
\] 
and
\[
\widetilde{U}_{b} = \widetilde{U}_{b}(X) := \{f \in X\kst (f\cup
b)\setminus \{v\} \not\in X \text{ for some } v\in b\}\,.
\]
Notice that
$U_b=\widetilde{U}_{b}=X$ unless $\partial b$ is a subcomplex of $X$.
Moreover, if $\partial b\subseteq X$, then
with $L_{b}:=\bigcap_{b^\prime\subset b}\link(b^\prime,X)$ we have
\begin{equation*}
X\setminus U_b = \begin{cases}  \emptyset \\\overline{\st}(b)
\end{cases}
\hspace{0.0em}\mbox{and}\hspace{0.7em}
X\setminus \widetilde{U}_{b} = \begin{cases}
\partial b \ast L_b& \text{if $b$ is a non-face},\\
(\partial b \ast L_b) \cup \overline{\st}(b)& \text{if $b$ is a face}.
\end{cases}
\vspace{1ex}
\end{equation*}

{{\bf Example \ref{oct}.7} (continued) }
Going back to the degree $\emptyset-(yy')$ from
Example \ref{oct}.5,
we see that $\XD\setminus U_b$ is the (closed) subcomplex induced from the edge
$(yy')$. On the other hand, $L_b$ is the boundary of the rectangle
$(xzx'z')$, hence $\partial b \ast L_b$ is the octahedron (without the
diagonals). Thus, while $U_b$ is $\XD$ with the closed diagonal $yy'$
being removed, $\widetilde{U}_b$ consists of the (disjoint) union of the
open diagonals $(xx')$ and $(zz')$. Comparing this with the easy
$(N,\widetilde{N})$ in Example \ref{oct}.5,
we feel that one has to pay for the advantage of getting open sets. 
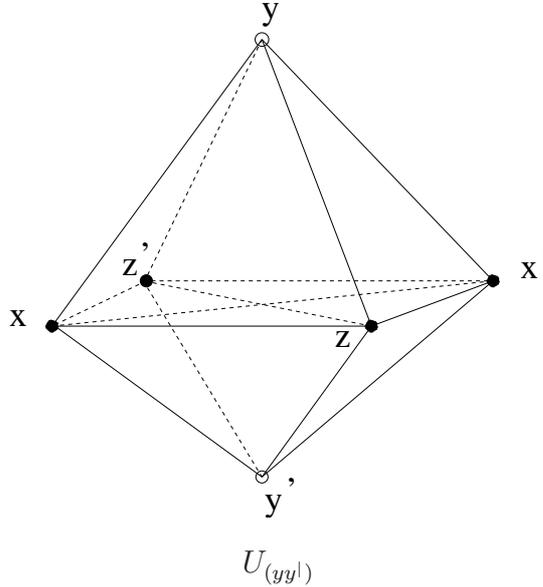
\vspace{1ex}

$$\begin{picture}(0,0)%
\includegraphics{diamant2.pstex}%
\end{picture}%
\setlength{\unitlength}{1657sp}%
\begingroup\makeatletter\ifx\SetFigFont\undefined
\def\x#1#2#3#4#5#6#7\relax{\def\x{#1#2#3#4#5#6}}%
\expandafter\x\fmtname xxxxxx\relax \def\y{splain}%
\ifx\x\y   
\gdef\SetFigFont#1#2#3{%
  \ifnum #1<17\tiny\else \ifnum #1<20\small\else
  \ifnum #1<24\normalsize\else \ifnum #1<29\large\else
  \ifnum #1<34\Large\else \ifnum #1<41\LARGE\else
     \huge\fi\fi\fi\fi\fi\fi
  \csname #3\endcsname}%
\else
\gdef\SetFigFont#1#2#3{\begingroup
  \count@#1\relax \ifnum 25<\count@\count@25\fi
  \def\x{\endgroup\@setsize\SetFigFont{#2pt}}%
  \expandafter\x
    \csname \romannumeral\the\count@ pt\expandafter\endcsname
    \csname @\romannumeral\the\count@ pt\endcsname
  \csname #3\endcsname}%
\fi
\fi\endgroup
\begin{picture}(7967,7731)(2296,-7414)
\end{picture}
$$
$$U_{(yy^|)}$$

\par

\begin{lemma}\klabel{opensets} 
$\;H^{\bullet}(\langle N_{\emptyset -
b}\rangle,\langle\widetilde{N}_{\emptyset - b}\rangle,\kk) 
\;\simeq\;
H^{\bullet}(\langle U_{b}\rangle, \langle \widetilde{U}_{b}\rangle,\kk)
$.
\end{lemma}

\begin{proof} If $f$ is in $U$ (respectively $\widetilde{U}$), then
$f\setminus b$ is in $N$ (respectively $\widetilde{N}$).  Now, if
$f\setminus b\ne \emptyset$ for all $f\in U$, then there is a standard
retraction taking $\alpha \in \langle f \rangle$ to an $\alpha^\prime
\in \langle f\setminus b \rangle$ which fits together to make
$(\langle N\rangle,\langle\widetilde{N}\rangle)$ a strong deformation
retract of $(\langle U\rangle,\langle\widetilde{U}\rangle)$.  (See
e.g.\  \cite[Proof of 3.3.11]{spa:alg}.)\\
If $b$ is a face, then $f\setminus b=\emptyset$ is
impossible, since then $f\cup b=b\in X$.  If $b$ is a non-face and
$\emptyset \in \widetilde{N}$, i.e.,  $\partial b \not \subseteq X$,
then all four spaces are cones and there is nothing to prove.  If $b$
is a non-face and $\emptyset \not\in \widetilde{N}$, then both
$\langle N \rangle$ and $\langle U \rangle$ are cones, so
$H^{i}(\langle N \rangle,\langle\widetilde{N}\rangle) \simeq
\widetilde{H}^{i-1}(\langle\widetilde{N}\rangle)$ and $H^{i}(\langle U
\rangle,\langle\widetilde{U}\rangle) \simeq
\widetilde{H}^{i-1}(\langle\widetilde{U}\rangle)$ and the above
retraction works again.
\end{proof}

When we plug this into Theorem~\ref{HNN} and use
Proposition~\ref{aempty} we get as a corollary the following 
description of the graded pieces of $T^i_{A_X}$.

\begin{theorem}\klabel{topopen}
The homogeneous pieces in degree $\mathbf{c}=\mathbf{a}-\mathbf{b}$
(with disjoint supports $a$ and $b$) of the cotangent cohomology of the
Stanley-Reisner ring $A_X$ vanish unless $a\in X$ 
and $\mathbf{b} =b\ne \emptyset$.  If these conditions are satisfied,
we have isomorphisms 
\[
T^i_\kc(X)
\;\simeq\; 
H^{i-1}\big(\langle U_{b}(\link(a,X))\rangle, \,
\langle \widetilde{U}_{b}(\link(a,X))\rangle,\,\kk\big)
\;\text{ for } i=1,2
\]
unless $b$ consists of a single vertex.  If $b$ consists of only one
vertex, then the above formulae become true if we use the reduced
cohomology instead.
\end{theorem}
\par

The reduction to the $a=\emptyset$ case 
also appears
in a completely different context.
One of the main issues of the paper \cite{ac:def} is the 
deformation theory of $\PP(X)\subseteq\PP^n$.
The relation to the deformation theory of its affine charts
$D_{+}(x_v)$ is governed by the localization maps
$T^i_{A_X,0}\to T^i_{D_{+}(x_v)}$; here ``$0$'' is meant with respect to the
usual $\Z$-grading of $A_X$, 
%
%
and the localization maps are obtained by dehomogenizing.
Now, the point is that these affine charts also come from simplicial complexes.
If $v\in[n]$, then $D_{+}(x_v)=\A(\link(v,X))$,
and we can use the techniques developed so far to describe the 
localization maps.
\par
 
\begin{remark}
Although, strictly speaking, we should consider $\link(v,X)$ as a
subcomplex of $\Delta_n$, 
the $T^i$ depend only on the
complex itself. In particular, when we look at graded parts, 
Proposition \ref{aempty}
shows that $T^i_{\ka-b}(\link(v))=0$ if $a$ or $b$ contain
non-vertices of $\link(v,X)$.
\end{remark}

\begin{lemma}\klabel{local} 
Let $\kc=\ka-b$ belong to degree $0$, i.e.,
$\deg \ka=\#b$ with $\deg$ denoting the sum of entries.
Fix a vertex $v\in[n]$.
\vspace{0.5ex}\\
{\rm (i)}
Localization with respect to $v$ maps 
$\,T^i_{\ka-b}\subseteq T^i_{A_X,0}$ into the graded summand
$\,T^i_{(\ka\setminus v)-(b\setminus v)}(\link(v,X))
\subseteq T^i_{D_{+}(x_v)}$
where $(\ka\setminus v)$ means cancellation of the $v$ entry.
\vspace{0.5ex}\\
{\rm (ii)}
The map of (i) is induced from
$\,\psi_{a-b}(v):
N_{(a\setminus v)-(b\setminus v)}(\link(v,X)) \rightarrow N_{a-b}(X)$
defined by
$\psi_{a-b}(v)(g) = g \cup (v \setminus b)$. 
It is compatible with the $\widetilde{N}$ level.
\vspace{0.5ex}\\
{\rm (iii)}
If $v\in a$, then the localization map is an isomorphism 
in degree $\ka-b$.
\end{lemma}

\begin{proof} It is straightforward to check that $\psi_{a-b}(v)$ is
well defined and has the necessary properties to induce
$\,\psi^\ast_{a-b}(v):T^i_{\ka-b}(X)\to 
T^i_{(\ka\setminus v)-(b\setminus v)}(\link(v,X))\,$ by Theorem~\ref{HNN}.
Moreover, it is clear that this means exactly
dehomogenization with respect to $v$, hence 
$\psi^\ast_{a-b}(v)$ coincides with the localization map.\\
Finally, to prove that 
$\,\psi_{a-b}(v): g\mapsto g\cup v$ is an isomorphism in case of $v\in a$, 
we use that 
$\,\link\big((a\setminus v),\, \link(v,X)\big) = \link(a,X)$
and apply Proposition \ref{aempty} to both 
$N_{(a\setminus v)-b}(\link(v,X))$ and $N_{a-b}(X)$.
%
\end{proof}

In fact, localizing with respect to {\em all} variables $x_v$ with $v$ 
running through the vertices of a given face $a\in X$ is induced
by the map
$\,\psi_{a-b}(a):
N_{\emptyset-b}(\link(a,X)) \rightarrow N_{a-b}(X)$
sending $g$ to $g \cup a$. This is the inverse of the $a$ killing map of
Proposition~\ref{aempty}.
%


\begin{theorem}\klabel{inject} 
The maps 
$\,T^i_{A_X,0}\to \bigoplus_{v\in [n]} T^i_{D_{+}(x_v)}$
are injective for $i=1,2$.
\end{theorem}

\begin{proof} 
If a graded piece $T^i_{\ka'-b'}(\link(v))$ meets
the image of $T^i_0(X) \rightarrow T^i(\link(v))$, then 
its pre-image is a unique summand
$T^i_{\ka-b}(X)\subseteq T^i_0(X)$. 
Indeed, by Proposition~\ref{local}, the only possibility to
get degree $0$ is
$\,\ka=\ka^\prime + (\#b^\prime-\deg\ka^\prime) 
\,v$,
$\,b=b'$ if $\deg \ka'\leq \#b'$
and 
$\ka=\ka'$, $\,b=b'\cup v$ if $\deg \ka'>\#b'$.
This means that images of elements with different multidegree
$\kc$ cannot cancel each other, and it remains to consider the
multigraded pieces
\vspace{-0.5ex}
\[
\,\mbox{$\bigoplus_{v\in [n]}$} \psi^\ast_{a-b}(v):\;
T^i_{\ka-b}(X) \longrightarrow 
\mbox{$\bigoplus_{v\in [n]}$} 
T^i_{(\ka\setminus v)-(b\setminus v)}(\link(v))
\]
with $\deg\ka = \#b$.  
Since, by Proposition~\ref{local}(iii), every summand
$\,\psi^\ast_{a-b}(v)$ with $v\in a$ is an isomorphism,
we obtain the injectivity of the above map whenever
$a$ has vertices at all.
On the other hand, since $\deg\ka = \#b$, the face $a$ cannot be
empty.
\end{proof}


\section{Examples}
\klabel{ex}

First, in Examples \ref{S0} and \ref{En}, we present 
the complete treatment of the easiest $X$ of all, the triangulations of  
0- and 1-dimensional manifolds. While the dimensions of $T^i$ for, say, the cone
over the $n$-gon are already well known, we can demonstrate how the 
multigrading comes in. Moreover, for higher-dimensional examples like 
the surfaces in Example \ref{dim2}, the smaller ones are needed because
they occur as links.


\begin{example}\klabel{S0}
Let $S^0=\{\emptyset,0,1\}$ be the $0$-dimensional complex 
consisting of two points only.
It may be considered a triangulation of the $0$-dimensional sphere.\\
If $b=\{0\}$, then $\widetilde{N}_{a-0}=\emptyset$ and 
$N_{a-0}=\{f\in S^0\kst a\subseteq f,\; f\cup 0\notin S^0\}=\{1\}$ 
for both possibilities $a=\emptyset$ or $a=\{1\}$.
In particular, $T^1_{a-0}(S^0)=T^2_{a-0}(S^0)=0$.\\
If $b=\{0,1\}$, then $a=\emptyset$, hence 
$N_{\emptyset-\{0,1\}}=\{\emptyset\}$ and 
$\widetilde{N}_{\emptyset-\{0,1\}}=\emptyset$.
This yields $T^2_{\emptyset-\{0,1\}}(S^0)=0$, but 
\fbox{$T^1_{\emptyset-\{0,1\}}(S^0)$ is one-dimensional}.
\vspace{0.5ex}\\
How does this infinitesimal deformation perturb the $S^0$-equation
$x_0x_1=0$?
If $\varepsilon$ denotes the infinitesimal parameter from
$\C[\varepsilon]/ \varepsilon^2$, then one obtains
$x_0x_1-\varepsilon=0$.
\end{example}


\begin{example}\klabel{En}
Denote by $E_n$ the simplicial complex representing an $n$-gon
with $n\geq 3$.
Index the vertices cyclically with
$0,\dots, n-1$; all addition is done modulo $n$.\\
First, we will show how to use the $a=\emptyset$ reduction
from Proposition \ref{aempty}.
If $a$ is an edge, then $\link(a,E_n)=\emptyset$, hence 
$T^i_{\ka-b}=0$. 
If $a=\{1\}$ is a vertex, then 
$\link(a,E_n)=\{\emptyset,0,2\}\cong S^0$,
hence from Example~\ref{S0} we obtain 
\fbox{$\dim T^1_{\{1\}-\{0,2\}}(E_n)=1$} as the only non-trivial contribution;
it translates into $\,x_0x_2-\varepsilon x_1^{\geq 1}$.
\\
Let us now assume that $a=\emptyset$.
If $\#b=1$, i.e., $b$ is a vertex, 
then $\widetilde{N}_{\emptyset-b}=\emptyset$,
and $\langle N_{\emptyset-b}\rangle$ equals $|E_n|$ after removing the edges
containing $b$. In particular, $\langle N_{\emptyset-b}\rangle$ is
contractible,
and the corresponding $T^i_{\emptyset-b}(E_n)$ are trivial.
If $b$ is an edge, then $\langle N_{\emptyset-b}\rangle$ looks similar,
and $\langle\widetilde{N}_{\emptyset-b}\rangle$ equals
$\langle N_{\emptyset-b}\rangle$ without the endpoints. 
We obtain $T^i_{\emptyset-b}(E_n)=0$ for $n\geq 4$, but
\fbox{$\dim T^1_{\emptyset-b}(E_3)=1$} yielding $x_0x_1x_2-\varepsilon x_i$.
\\
The final case is that $a=\emptyset$ and $b\notin E_n$.
If $\#b\geq 3$, then, except for $n=3$, the set
$b$ always contains proper subsets which are non-faces. In particular,
$\widetilde{N}_{\emptyset-b}=N_{\emptyset-b}$, leading to 
$T^i_{\emptyset-b}(E_n)=0$. The exception is
\fbox{$\dim T^1_{\emptyset-\{0,1,2\}}(E_3)=1$}
yielding $x_0x_1x_2-\varepsilon$.
It remains to take two non-adjacent vertices for $b$, say $b=\{u,v\}$.
Since $\emptyset\in N_{\emptyset-b}$, the set $\langle N_{\emptyset-b}\rangle$ 
is always a cone, hence contractible. 
In particular, $T^1_{\emptyset-b}(E_n)=0$ if 
$\widetilde{N}_{\emptyset-b}\neq\emptyset$, and this is always the case 
except for \fbox{$\dim T^1_{\emptyset-\{v,v+2\}}(E_4)=1$};
in terms of equations: $\,x_vx_{v+2}-\varepsilon$.
On the other hand, the long exact cohomology sequence for the pair
$(N_{\emptyset-b}, \widetilde{N}_{\emptyset-b})$ 
yields $\,T^2_{\emptyset-b}(E_n)=
\widetilde{H}^0(\widetilde{N}_{\emptyset-b})$.
Since the set $\langle\widetilde{N}_{\emptyset-b}\rangle$ equals $|E_n|$ with
$u$ and $v$ and the adjacent edges being removed,
we eventually obtain
\fbox{$\dim T^2_{\emptyset-\{u,v\}}(E_n)=1$ whenever $|u-v|\geq 3$}. 
\vspace{1ex}\\
Adding up, we find that $T^2(E_n)=0$ if $n\le 5$, and that 
$\dim T^2(E_n)=n(n-5)/2$ if $n\ge 6$.
In the latter case, we can even locate where the cup product takes place.
Considering the coarse $\Z$-grading, we see that 
$T^1(E_n)$ spreads in degre $\geq -1$, and $T^2(E_n)$ sits in degree $-2$.
Hence, the cup product lives in the pieces
$T^1_{-1}\times T^1_{-1}\to T^2_{-2}$ only.
Using the $\Z^n$ multigrading, one obtains a finer result.
The cup product splits into products
\vspace{-1ex}
\[
T^1_{\{v\}-\{v-1,v+1\}}(E_n)\;\times\; T^1_{\{v+1\}-\{v,v+2\}}(E_n)
\;\longrightarrow\; T^2_{\emptyset-\{v-1,v+2\}}(E_n)
\]
with all the three vector spaces being one-dimensional.
See the continuation on p.~\pageref{excup} for more information.
\end{example}


\begin{example}\klabel{dim2}
Let us look at the degree $0$ deformations when $X$ is 
the triangulation of a two-dimensional manifold.
The non-zero multigraded pieces $T^1_{\ka-b}(X)$ of $T^1_0(X)$
require $a\neq\emptyset$, hence they are induced from lower-dimensional links.
They are all $1$-dimensional
and are gathered in the following list.
\vspace{-1ex}
\par

\begin{list}{\textup{(\roman{temp})}}{\usecounter{temp}}
\item\parbox[t]{7cm}{
$a$ is an edge; $|a|=|b|=2$.\\
${\mathbf a}=a$.
}
{
\unitlength=0.5pt
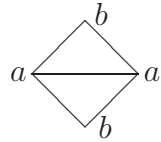
\begin{picture}(0.00,0.00)(-30.00,100.00)
\put(52.00,125.00){\makebox(0.00,0.00){$b$}}
\put(55.00,40.00){\makebox(0.00,0.00){$b$}}
\put(90.00,80.00){\makebox(0.00,0.00){$a$}}
\put(-10.00,80.00){\makebox(0.00,0.00){$a$}}
\put(0.00,80.00){\line(1,0){80.00}}
\put(40.00,40.00){\line(-1,1){40.00}}
\put(80.00,80.00){\line(-1,-1){40.00}}
\put(40.00,120.00){\line(1,-1){40.00}}
\put(0.00,80.00){\line(1,1){40.00}}
\end{picture}}

\item
\parbox[t]{7cm}{
$a$ is a vertex with valency $4$; \\ $|a|=1,\,|b|=2$.\\
${\mathbf a}=2\cdot a$.
}
{
\unitlength=0.5pt
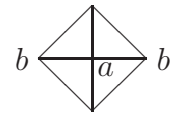
\begin{picture}(0.00,0.00)(-160.00,80.00)
\put(50.00,52.00){\makebox(0.00,0.00){$a$}}
\put(50.00,5.00){\makebox(0.00,0.00){}}
\put(50.00,110.00){\makebox(0.00,0.00){}}
\put(93.00,60.00){\makebox(0.00,0.00){$b$}}
\put(-13.00,60.00){\makebox(0.00,0.00){$b$}}
\put(40.00,60.00){\line(0,-1){40.00}}
\put(40.00,100.00){\line(0,-1){40.00}}
\put(40.00,60.00){\line(1,0){40.00}}
\put(0.00,60.00){\line(1,0){40.00}}
\put(40.00,20.00){\line(-1,1){40.00}}
\put(80.00,60.00){\line(-1,-1){40.00}}
\put(40.00,100.00){\line(1,-1){40.00}}
\put(0.00,60.00){\line(1,1){40.00}}
\end{picture}}

\item
\parbox[t]{7cm}{
$a$ is a vertex with valency $3$;\\ $|a|=1,\,|b|=3$;\\
${\mathbf a}=3\cdot a$.
}
{
\unitlength=0.5pt
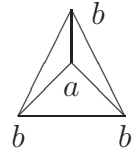
\begin{picture}(0.00,0.00)(-30.00,80.00)
\put(60.00,98.00){\makebox(0.00,0.00){$b$}}
\put(40.00,40.00){\makebox(0.00,0.00){$a$}}
\put(80.00,5.00){\makebox(0.00,0.00){$b$}}
\put(0.00,5.00){\makebox(0.00,0.00){$b$}}
\put(40.00,60.00){\line(0,1){40.00}}
\put(40.00,100.00){\line(1,-2){40.00}}
\put(0.00,20.00){\line(1,2){40.00}}
\put(80.00,20.00){\line(-1,0){80.00}}
\put(40.00,60.00){\line(1,-1){40.00}}
\put(0.00,20.00){\line(1,1){40.00}}
\end{picture}}

\item
\parbox[t]{7cm}{
$a$ is a vertex with valency $3$;\\ $|a|=1,\,|b|=2$.\\
${\mathbf a}=2\cdot a$.
}
{
\unitlength=0.5pt
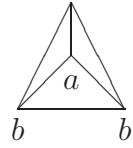
\begin{picture}(0.00,0.00)(-160.00,80.00)
\put(60.00,98.00){\makebox(0.00,0.00){}}
\put(40.00,40.00){\makebox(0.00,0.00){$a$}}
\put(80.00,5.00){\makebox(0.00,0.00){$b$}}
\put(0.00,5.00){\makebox(0.00,0.00){$b$}}
\put(40.00,60.00){\line(0,1){40.00}}
\put(40.00,100.00){\line(1,-2){40.00}}
\put(0.00,20.00){\line(1,2){40.00}}
\put(80.00,20.00){\line(-1,0){80.00}}
\put(40.00,60.00){\line(1,-1){40.00}}
\put(0.00,20.00){\line(1,1){40.00}}
\end{picture}}

\end{list}

The perturbation of the equations also comes from the
corresponding one of the lower-dimensional link.
E.g., denote the vertices in (i) such that $a=\{y,x_1\}$
and $b=\{x_0,x_2\}$. Then we obtain $x_0x_2-\varepsilon x_1y$. 

Moreover, $T^2_0(X)$ is only present for $a$ being a vertex of valency 
at least $6$, $\ka=2\cdot a$, and $b$ consisting of two vertices 
having exactly $a$ as a common neighbor.
\end{example}


{{\bf Example \ref{oct}.8} (finished) }
Finally, we conclude our running Example \ref{oct};
the simplicial complex $\XD$ is the octahedron plus the three diagonals
$(xx')$, $(yy')$, and $(zz')$.%
\vspace{0.3ex}\\
First, if $a$ is a non-empty face,
then one has the following types of links:
$\,\link(xyz)=\emptyset$, $\,\link(xx')=\emptyset$
(both yielding $T^1=T^2=0$),
$\,\link(xy)=\{z,z',\emptyset\}\cong S^0$ yielding
\fbox{$\dim T^1_{(xy)-(zz')}(\XD)=1$}.
Moreover, 
we studied the case $a=\{x\}$
in Example~\ref{oct}.6; the non-zero results were
\fbox{$\dim T^2_{(x)-(x')}(\XD)=1$}
and,
for any $\ast\in\{y^{(')}, z^{(')}\}$,
\fbox{$\dim T^1_{(x)-(\ast)}(\XD)=1$}.
In terms of equations, the $T^1$-contributions look like
$yzz'-\varepsilon xy^2$ in the first case,
and like $\,xx'y-\varepsilon x^2x'$, $\,xyy'-\varepsilon x^2y'$,
$\,yy'z-\varepsilon xy'z$ in degree $(x)-(y)$.
\vspace{0.3ex}\\
Second, if $a=\emptyset$, then there are two major cases to distinguish.
If $b$ is not a face, then one knows in general that
$\emptyset\in N_{\emptyset-b}(X)$, i.e., $\langle N_{\emptyset-b}(X)\rangle$
is contractible. Hence, by the long exact cohomology sequence
for the pair $(N,\widetilde{N})$,
$\dim T^1_\kc=1$ and $T^2_\kc=0$ if
$\widetilde{N}_{\emptyset-b}(X)=\emptyset$,
and, otherwise,
$T^1_\kc=0$, $\,T^2_\kc=\widetilde{H}^0(\widetilde{N}_{\emptyset-b}(X))$.
 In the case of $X=\XD$, this always yields $T^i_\kc=0$.
\vspace{0.3ex}\\
It remains to consider faces for $b$.
While $b=\{x'\}$ or $\{x',y'\}$ do not yield any $T^i$,
we obtained in Example~\ref{oct}.5 that
\fbox{$\dim T^2_{\emptyset-(yy')}(\XD)=4$}.
\vspace{1ex}
\par


\section{Appendix: The cup product}
\klabel{cup}

The cup product $T^1\times T^1 \rightarrow T^2$ is an important tool 
to obtain more information about deformation theory 
than just the knowledge of the tangent or obstruction spaces
$T^i$ themselves.
The associated quadratic form $T^1\to T^2$ 
describes the equations of the versal base space up to second order.
\par

In the case of Stanley-Reisner rings $A_X$, we only managed 
to get a nasty description of this product using the language of 
Proposition~\ref{TiN}.
We have not yet found a relation to the geometry of the complex.
However, since the cup product provides 
important information needed in some applications in \cite{ac:def},
we have decided to present it in an appendix.
It is suggested that
overly-sensitive readers quit reading at this point.
\par
 
The cup product can be defined in the
following way (see \cite[5.1.5]{la:for}):
Let $A=P/I$ with $I$ generated by equations $f^p$.  
If $\varphi\in \Hom(I, A)$, lift the images of the $f^p$ 
obtaining elements
$\widetilde{\varphi}(f^p)\in P$.  Given a relation $r\in R$, the
linear combination $\langle r, \widetilde{\varphi} \rangle :=
\sum_p r_p\,\widetilde{\varphi}(f^p)$ vanishes in $A$, i.e.\ it is
contained in $I$.  If $\varphi,\psi \in \Hom(I, A)$ represent two
elements of $T^1$, then we define for each relation $r\in R$
\begin{equation}\label{gen}
(\varphi \cup \psi)(r) := \psi(\langle r, \widetilde{\varphi}
\rangle) + \varphi(\langle r, \widetilde{\psi} \rangle)\, .
\end{equation}
This determines a well defined element of $T^2$.
If $A=A_X$ is a Stanley-Reisner ring, then
the cup product respects the multigrading, and,
using Proposition~\ref{TiN}, we can give a formula
for $\cup : T^1_{\mathbf{a}^1-\mathbf{b}^1} \times
T^1_{\mathbf{a}^2-\mathbf{b}^2} \to T^2_{\mathbf{a}-\mathbf{b}}$
with $\mathbf{a}-\mathbf{b} = \mathbf{a}^1-\mathbf{b}^1 +
\mathbf{a}^2-\mathbf{b}^2$.
\vspace{1ex}
\par


\begin{proposition}\klabel{cp} 
$\,\cup : T^1_{\mathbf{a}^1-\mathbf{b}^1} \times
T^1_{\mathbf{a}^2-\mathbf{b}^2} \to T^2_{\mathbf{a}-\mathbf{b}}$ is
determined by the following:
\vspace{0.5ex}\\
%
{\rm (i)}
If $b^1 \cap b^2 \ne \emptyset$, then $\cup = 0$.
\vspace{0.5ex}\\
%
{\rm (ii)}
If $b^1 \cap b^2 = \emptyset$, then
$b=(b^1\setminus a^2) \cup (b^2\setminus a^1)$
and
$a=(a^1\setminus b^2) \cup (a^2\setminus b^1)
\cup(a^1_{\geq 2}\cup a^2_{\geq 2})$
with $a^{\kbb}_{\geq 2}:=\{v\in [n]\,|\; a^{\kbb}_v \geq 2\}$ denoting
the locus of higher multiplicities.
\vspace{0.5ex}\\
%
{\rm (iii)}
If $(f,g) \in N_{a-b}^{(2)}$, choose maximal subsets $d,e
\subseteq b$ such that $f\cup (b\setminus d)$ and $g \cup (b\setminus
e) \notin X$. 
%
%
If $\varphi \in T^1_{a^1-b^1}$ and $\psi \in
T^1_{a^2-b^2}$, then the value of $(\varphi \cup \psi)(f,g)$ is
\vspace{-1ex}
\begin{multline*}
\Big[\varphi\Big([f\setminus b^1]\cup[b^2\setminus d]\Big)
- \varphi\Big([g\setminus b^1]\cup[b^2\setminus e]\Big)\Big]
\cdot \psi\Big(f\cup g \cup a^2\setminus b^2\Big) \\
  + \Big[\psi\Big([f\setminus b^2]\cup[b^1\setminus d]\Big)
- \psi\Big([g\setminus b^2]\cup[b^1\setminus e]\Big)\Big]
\cdot \varphi\Big(f\cup g \cup a^1\setminus b^1\Big)
\end{multline*}
with $\varphi, \psi$ defined to be zero on non-elements of $X$.
\end{proposition}

In fact, the maximality of $d$ and $e$ is not quite necessary. 
The point is to choose them non-empty whenever possible.
\par


\begin{proof}
(i) We know that
$\mathbf{b}\in\{0,1\}^n$ if $T^2_{\mathbf{a}-\mathbf{b}} \ne 0$. This
would not be the case if $b^1 \cap b^2 \ne \emptyset$. Statement (ii)
is a straightforward calculation.\\
For (iii), we must first recall what we did in Step 3 of the proof of
Proposition~\ref{TiN}. Antisymmetric functions $\lambdaM(\kbb,\kbb)$
on the $M$ level had been turned into antisymmetric
functions $\lambdaN(\kbb,\kbb)$ on the $N$ level
via $\lambdaN(f,g):=\lambdaM(m_f,m_g)$ with certain elements 
$m_{\kbb}\in \Phi^{-1}(\kbb)$.
Hence, setting
$p := m_f= f\cup (b\setminus d)$ and $q := m_g= g \cup (b\setminus e)$,
we may compute the value of 
$(\varphi \cup \psi)(f,g)$ 
by applying the expression~\ref{gen} on the relation $R_{p,q}$
described in the beginning of Section \ref{TComplex}.
Using
$\widetilde{\varphi}(\mathbf{x}^{p}) = 
\varphi(\Phi_{a^1-b^1}(p))\cdot\kx^{p+\ka^1-b^1}$,
we obtain
\[
\langle R_{p,q}, \widetilde{\varphi} \rangle
=\big[\varphi \big(\Phi_{a^1-b^1}(p)\big)
-\varphi \big(\Phi_{a^1-b^1}(q)\big)\big] \cdot 
\kx^{(p\cup q)+\ka^1-b^1}.
\]
Plugging this into \ref{gen}, we get
\begin{equation*}\begin{split}
(\varphi \cup \psi)(f,g) \;=\;&
\big[\varphi\big(\Phi_{a^1-b^1}(p)\big)-
\varphi\big(\Phi_{a^1-b^1}(q)\big)\big]\cdot
\psi\big(\Phi_{a^2-b^2}
\big[\Phi_{a^1-b^1}(p \cup q)\big]\big)\\
  \;+\;&
\big[\psi\big(\Phi_{a^2-b^2}(p)\big)-
\psi\big(\Phi_{a^2-b^2}(q)\big)\big]\cdot
\varphi\big(\Phi_{a^1-b^1}
\big[\Phi_{a^2-b^2}(p \cup q)\big]\big)\,.
\end{split}\end{equation*}
To finish the proof, it is still necessary to take a closer look at the
occurring arguments, i.e., to calculate
\[
\renewcommand{\arraystretch}{1.2}
\begin{array}{rl}
\Phi_{a^1-b^1}(p)
&=\Phi_{a^1-b^1}\big(f\cup(b\setminus d)\big)\\
&=
\big[f\cup(b^1\setminus a^2\setminus d) \cup
(b^2\setminus a^1\setminus d)
\cup a^1\big]\setminus b^1\\
&=
\big[f \cup(b^2\setminus d)\big]\setminus b^1\quad
(\text{since } (a^1\setminus b^2)\subseteq f
\text{ and } a^1\cap d=\emptyset)\\
&=
(f\setminus b^1)\cup(b^2\setminus d)
\end{array}
\vspace{-2ex}
\]
and
\begin{align*}
\Phi_{a^2-b^2}\big[\Phi_{a^1-b^1}(p \cup q)\big]&=
\Phi_{a^2-b^2}\big[\Phi_{a^1-b^1}\big((f\cup g) \cup(b\setminus
(d\cap e)\big)\big]\\
&=
\Phi_{a^2-b^2}\big[((f\cup g)\setminus b^1)\cup(b^2\setminus (d\cap e))\big]\\
&=
\big[((f\cup g\cup a^2)\setminus (b^1\setminus a^2)\big]
\setminus b^2\\
&= \big[f\cup g \cup a^2\big]\setminus b^2 \quad (\text{since }
(b^1\setminus a^2)\cap (f\cup g) =\emptyset)\, .
\end{align*}
One can check that all the arguments of $\varphi$ and $\psi$ are in
$N_{a^i-b^i}$ if they are in $X$.
\end{proof}


\klabel{excup}
{\bf Example \ref{En} (continued). }
We are going to calculate the cup product mentioned at the end of
Example~\ref{En}. 
With $f_1:=\{v-2\}$,  $f_2:=\{v\}$ we choose representatives from
both connected components of $\widetilde{N}_{\emptyset - \{v-1,v+2\}}$.
Since 
$N_{\emptyset - \{v-1,v+2\}} = 
\widetilde{N}_{\emptyset - \{v-1,v+2\}} \cup \{\emptyset\}$,
we obtain that $\lambdaN(f_1,\emptyset)$ and $\lambdaN(f_2,\emptyset)$
suffice to know about a function $\lambdaN:N_{\emptyset - \{v-1,v+2\}}\to\kk$
from Proposition~\ref{TiN}. Since $T^2$ results from dividing out a subspace,
we obtain $[\lambdaN(f_1,\emptyset)-\lambdaN(f_2,\emptyset)]$ as the ultimate
coordinate of $T^2_{\emptyset - \{v-1,v+2\}}(E_n)$.
As auxillary elements $d,e\subseteq b$ (cf.~Proposition~\ref{cp}),
we may choose $d:=\{v-1\}$ for both $f_1$ and $f_2$ and
$e:=\emptyset$ for the second argument $g:=\emptyset$.\\
On the other hand, we know that
$N_{\{v\}-\{v-1,v+1\}}=\{v\}$ and $N_{\{v+1\}-\{v,v+2\}}=\{v+1\}$
with $\widetilde{N}=\emptyset$ in both cases.
Hence, the $T^1$ spaces are represented by maps $\varphi$ and $\psi$
yielding $1$ on the faces $\{v\}$ and $\{v+1\}$, respectively.
Applying Proposition~\ref{cp}, we find
\[
\renewcommand{\arraystretch}{1.2}
\begin{array}{@{}l@{}l@{}}
(\varphi \cup \psi)(f_1,\emptyset) & \;=  \, 
[\varphi(\{v-2,v,v+2\}) - \varphi(\{v, v+2\})] \cdot 
\psi(\{v-2,v+2\}) \\
& \hspace*{1.5em}+ \, [\psi(\{v-2,v+2\}) - \psi(\{v-1,v+1\})] \cdot
\varphi(\{v-2,v\}) \;=\; 0
\end{array}
\vspace{-1ex}
\]
and
\[
\renewcommand{\arraystretch}{1.2}
\begin{array}{@{}l@{}l@{}}
(\varphi \cup \psi)(f_2,\emptyset) & \;=  \, 
[\varphi(\{v,v+2\}) - \varphi(\{v, v+2\})] \cdot 
\psi(\{v+1\}) \\
& \hspace*{1.5em}+ \, [\psi(\{v+1\}) - \psi(\{v-1,v+1\})] \cdot
\varphi(\{v\}) \;=\; 1\,.
\end{array}
\]
Thus, the cup product mentioned at the end of
Example~\ref{En} yields $\varphi\cup\psi=1$.
\par

\begin{corollary}\klabel{7gon}
Let $n\geq 7$. If $t_1,\dots,t_n$ denote the coordinates of 
$\,T^1_{-1}(E_n)=
\oplus_{v\in\Z/n\Z} \,T^1_{\{v\}-\{v-1,v+1\}}(E_n)$,
then the equations of the negative part of the base space $S$ of the
versal deformation of $E_n$ are $t_vt_{v+1}=0$ for $v\in\Z/n\Z$.
In particular, $E_n$ is not smoothable over $S$.
\end{corollary} 

\begin{proof} 
Via the cup product, we see that each part
$T^2_{\emptyset-\{v-1,v+2\}}(E_n)$ is responsible for $t_vt_{v+1}$
in the quadratic part of the obstruction equations.
Moreover,
since $T^2$ is concentrated in degree $-2$, no higher order
obstructions involving {\em only} degree $-1$ deformations can appear,
i.e., $S$ is described by the desired equations.\\
Thus, in any flat deformation of degree $-1$, every other
parameter must vanish. 
One directly checks that this means that any fiber is
singular, in fact reducible.
\end{proof}

In contrast, if $n=6$, then each of the three $T^2$-pieces is the 
common target of two different cup products. In particular,
the negative part of the base space $S$ is given by the equations
$t_0t_1-t_3t_4=t_1t_2-t_4t_5=t_2t_3-t_5t_0=0$. This yields the cone
over the three-dimensional, smooth, projective toric variety 
induced by the prism over the standard triangle.
\par


\bibliographystyle{amsalpha}
{\small
\providecommand{\bysame}{\leavevmode\hbox to3em{\hrulefill}\thinspace}

}

{\small
\setbox0\hbox{FB Mathematik und Informatik, WE2}
\parbox{\wd0}{
Klaus Altmann\\
FB Mathematik und Informatik, WE2\\
Freie Universit\"at Berlin\\
Arnimallee 3\\
D-14195 Berlin, Germany\\
email: altmann@math.fu-berlin.de}
\setbox1\hbox{University of Oslo at Blindern}\hfill\parbox{\wd1}{
Jan Arthur Christophersen\\
Department of Mathematics\\
University of Oslo at Blindern\\
Oslo, Norway\\
email: christop@math.uio.no}}

\end{document}